\documentclass{amsart}

\usepackage[latin1]{inputenc}
\usepackage{amsmath, amsthm, amssymb, booktabs, multirow, graphicx, booktabs}
\usepackage[colorlinks]{hyperref}
\usepackage{todonotes}
\usepackage{array}
\usepackage{subcaption}
\usepackage{enumitem}
\usepackage{float}

\newtheorem{defin}{Definition}[section]

\newtheorem{theorem}[defin]{Theorem}
\newtheorem{lemma}[defin]{Lemma}

\newcommand{\R}{\mathbb{R}}

\newcommand{\F}{\mathbb{F}}
\newcommand{\Q}{\mathbb{Q}}

\newcommand{\N}{\mathbb{N}}

\DeclareMathOperator{\Tr}{Tr}

\DeclareMathOperator{\Gram}{Gram}

\DeclareMathOperator{\Gal}{Gal}

\newcommand{\OptimalNumber}{{\mathbf N}^\star}
\newcommand{\newlowerbound}{t_{IS}}

\begin{document}

\title{Refinements on higher order Weil-Oesterl\'e bounds via a Serre type argument}

\address{E.~Hallouin, Institut de Mathe\'ematiques de Toulouse, UMR
  5219, UT2J, 31058 Toulouse, France}

\address{P.~Moustrou, Institut de Mathe\'ematiques de Toulouse, UMR
  5219, UT2J, 31058 Toulouse, France}

\address{M.~Perret, Institut de Mathe\'ematiques de Toulouse, UMR
  5219, UT2J, 31058 Toulouse, France}

\author{Emmanuel Hallouin}
\email{hallouin@univ-tlse2.fr}
  
\author{Philippe Moustrou}
\email{philippe.moustrou@math.univ-toulouse.fr}

\author{Marc Perret}
\email{perret@univ-tlse2.fr}

\date{\today}


\thanks{This work was supported by the French \textit{Agence
Nationale de la Recherche} project ANR-21-CE39-0009-BARRACUDA}

\begin{abstract}
Weil's theorem gives the most standard bound on the number of points of a curve over a finite field. 
This bound was improved by Ihara and Oesterl\'e for larger genus. 
Recently, Hallouin and Perret gave a new point of view on these bounds, that can be obtained by solving a sequence of semi-definite programs, and the two first steps of this hierarchy recover Weil's and Ihara's bounds. 
On the other hand, by taking into account arithmetic constraints, Serre obtained a refinement on Weil's bound. 
In this article, we combine these two approaches and propose a strengthening of Ihara's bound, based on an argument similar to Serre's refinement. 
We show that this generically improves upon Ihara's bound, even in the range where it was the best bound so far.
Finally we discuss possible extensions to higher order Weil-Oesterl\'e bounds. 
\end{abstract}

\maketitle

\markboth{E. Hallouin, P. Moustrou, and M. Perret}{Refinements on higher order Weil-Oesterl\'e bounds via a Serre type argument}

\setcounter{tocdepth}{1} 


\section*{Introduction}

This work deals with the problem of bounding above the number of rational points of an absolutely, irreducible, smooth, projective
curve~$X$, of genus~$g$, defined over the finite field~$\F_q$, where~$q$ is a prime power.
In the~$1940$s Weil \cite{Weil1,Weil2} proved that this
number~$\sharp X(\F_q)$ satisfies the inequalities
$$
(q+1) -2g\sqrt{q} \leq \sharp X(\F_q) \leq (q+1) + 2g\sqrt{q}.
$$
In~$1985$, Serre brought a new focus to this issue by giving a lecture on this topic at Harvard. The Gouvea's handwritten notes on these
courses, which have long circulated in the community, have recently been published \cite{SerreRP}. It contains lots of ideas to
improve the Weil bounds in several directions such as computing the exact values of the constants~$N_q(g)$ defined by
\[
N_q(g) = \max\left\{\sharp X(\F_q), \; \text{\begin{tabular}{l}$X$ absolutely, irreducible, smooth, projective\\$\F_q$-curve
of genus~$g$\end{tabular}}\right\}
\]
for small values of the genus~$g$, or computing better upper bounds for large genus. This course has opened the way for many new works
and developments leading to better and better bounds for different constants~$N_q(g)$
(on this topic or some other developments, see for example \cite{AubryIezzi,HoweF4,MR3349314,MR2961401,HL03,HL12,HL07,MR4876845} and the references in \cite{SerreRP}).
These improvements are now listed on the very useful website {\tt manYPoints} \cite{MP}, where the best-known lower and upper bounds
of the~$N_q(g)$'s for
small values of~$q$ and~$g$ are regularly updated.

This work is no exception to this rule, and continues the ideas developed in Serre's course, particularly two strategies
that we now point out.

First, taking into account the fact that the eigenvalues of the Frobenius are algebraic integers, Serre improved the Weil
bounds \cite{SerreWS} to what is now called the {\bf\em Weil-Serre bound}
(see section~\ref{sec:SerreWeil}):
$$
 (q+1) -g\lfloor2\sqrt{q}\rfloor \leq \sharp X(\F_q) \leq  (q+1) + g\lfloor2\sqrt{q}\rfloor.
$$
Of course this is only an improvement for non squares values of~$q$ and the gain is in~$g\{2\sqrt{q}\}$ where, for~$x\in\R$,
$\{x\}$ denotes its fractional part.

On the other hand, Ihara \cite{Ihara} noted that the Weil bounds cannot be optimal for large genus and Serre, in his course,
developed the so-called \emph{explicit formul\ae}.
This approach results in an optimization problem whose solution gives an upper bound for~$N_q(g)$, for increasing genus. Oesterl\'e solved this problem, leading to the Oesterl\'e bounds, see \cite[Chapter VI]{SerreRP}. 
In \cite{HP19}, the first and third authors gave a new point
of view on these bounds, by reproving them in the spirit of the original proof of the Weil bounds using intersection theory in the algebraic surface~$X\times X$. 
In their setting, recalled in section~\ref{sWn}, it is shown that there exists an explicit strictly increasing sequence~$(g_n)_{n\geq 1}$ of
non-negative real numbers and a sequence~$(\OptimalNumber_n)_{n\geq 1}$  of
strictly increasing functions from~$\left[g_n,+\infty\right)$ to~$\R$, such that for any~$g\geq g_n$, the value~$\OptimalNumber_n(g)$ is
an upper bound for~$N_q(g)$, that we call the {\bf\em Weil-Oesterl\'e bound} of order $n$. 
Moreover, the bound~$\OptimalNumber_{n+1}(g)$ is sharper than~$\OptimalNumber_n(g)$ for~$g > g_{n+1}$.
Furthermore, the bound~$\OptimalNumber_1(g)$ is nothing else than the Weil bound, while~$\OptimalNumber_2(g)$ is the Ihara bound \cite{Ihara}
and one has
\begin{align*}
&g_1=0,
&
&g_2 = \frac{\sqrt{q}\left(\sqrt{q}-1\right)}{2},
&
&g_3 = \frac{\sqrt{q}\left(q-1\right)}{\sqrt{2}}.
\end{align*}

The main contribution of this article is contained in section~\ref{sIS}. 
We combine the two previous strategies to obtain an improvement of Ihara's bound following Serre's improvement on Weil's bound. This leads us to a new upper bound for~$N_q(g)$ (Theorem~\ref{thm:main}) which works for every~$q$ and the gain compared to Ihara's bound is explicit.
Like Serre, this gain depends on the fractional part of a quantity
which is a little bit more difficult to analyse than for the Weil-Serre bound. 
Section~\ref{sEtudeGain} is devoted to this analysis, in several respects. First in Section~\ref{sec:qfix}, we fix~$q$ and let~$g$ go to infinity. We prove that, except for~$q=3$, our gain has an asymptote in~$g$ with positive sloape, and thus tends to infinity. 
Second we focus on the Ihara range~$[g_2,g_3]$. 
It is fair to compare our bound with Ihara's therein, because~$\OptimalNumber_2(g)$ does not get improved by Weil-Oesterl\'e bounds of higher order. 
In this interval, the fractional part in our bound makes more difficult the analysis of the gain. 
However, we use numerical experiments, for several values of~$q$, to compare our bound with Ihara's (see Section~\ref{sec:CompIhara}), and provide an explicit infinite sequence of couples $(q,g_q)$ with $g_q\in [g_2,g_3]$ such that this gain is the best we can hope for, see Section~\ref{sec:CompSeq}.  
Last in Section~\ref{sec:MP} we compare our bound with the entries in {\tt manYPoints}. 
In the Ihara range, we find more than~$150$ couples $(q,g)$ for which our bound improves upon Ihara's by~$1$. 
Among them, we recover the current record in more than~$130$ cases (where Ihara's bound was already improved by other techniques, mainly from \cite{HL03, HL12}), and improve this record for~$20$ couples $(q,g)$. 

We end by a more prospective Section~\ref{sec:extensions}. It is natural to ask whether or not Weil-Oesterl\'e bounds of higher order
could be improved in the same way. We give some experiments for the Weil-Oesterl\'e bound of order~$3$ which lead to new records, but we discuss
why our method fails to be easily generalized in greater order.

\section{Weil-Oesterl\'e bounds and Serre's trick to improve Weil's bound}\label{sOS}

Let $X$ be an absolutely irreducible smooth projective curve of genus $g$ defined over $\F_q$. 
Weil proved the existence of algebraic integers $\omega_1, \ldots, \omega_g$ of modulus $\sqrt{q}$ such that the number $N_k$ of $\F_{q^k}$-points of $X$ satisfies
\[
N_k  = \left(q^k+1\right) - \sum_{j=1}^g \omega_j^k + \overline{\omega_j}^k.
\]
Moreover, the family $\{ \omega_1, \ldots, \omega_g, \overline{\omega_1}, \ldots, \overline{\omega_g} \}$ is stable under the action of $\Gal(\overline{\Q}/\Q)$, they are nothing else than the eigenvalues of the Frobenius.
This results immediately implies bounds on $N_k$. 
Indeed, writing $\omega_j = \sqrt{q}e^{i\theta_j}$, one gets
\[
\left| N_k - (1+q^k) \right| \leq 2q^{k/2}\sum_{j=1}^g \left| \cos(k\theta_j) \right| \leq 2gq^{k/2}. 
\]
In particular, for $k=1$, one gets the celebrated Weil inequality
\begin{equation}\label{eq:Weil1}
\left| N_1 - (1+q) \right| \leq 2g\sqrt{q}.
\end{equation}
More generally, if we introduce 
\begin{equation}\label{eq:deftk}
t_k = \sum_{j=1}^g \omega_j^k + \overline{\omega_j}^k = 1+q^k - N_k,
\end{equation}
any lower bound on $t_k$ gives an upper bound on $N_k$. 

\subsection{Weil-Oesterl\'e bounds} \label{sWn}

Using the method of \emph{explicit formul\ae} (see \cite[Chapter V.3]{SerreRP}), any trigonometric polynomial $f$ which is nonnegative on the unit circle and whose coefficients in the cosine expansion are nonnegative gives a lower bound on $t_1$. 
Finding the best $f$ possible then becomes a conic optimization problem, and Oesterl\'e gave an explicit procedure to solve it, see \cite[Chapter VI]{SerreRP}.
In \cite{HP19}, the first and the third authors provided the following new point of view on this question.
Thanks to the Hodge index Theorem on the surface $X \times X$, the opposite of the intersection pairing defines a scalar product on the orthogonal of the space generated by horizontal and vertical divisors inside the space of divisors of~$X\times X$ up to numerical equivalence.
In this Euclidean space, if we denote by~$p(\Gamma^i)$, $i\geq 0$ ($\Gamma^0$ is the diagonal of~$X\times X$) the orthogonal projection of $\Gamma^i$ onto the orthogonal of the vertical and horizontal
parts, then one shows using elementary intersection theory that
\begin{equation}\label{eq:Toep}
\Gram\left(p(\Gamma^0),\ldots,p(\Gamma^n)\right)
=
\begin{pmatrix}
2g & t_1 & t_2 & \ldots  & t_n \\
t_1 & 2gq & qt_1 & \ddots  & qt_{n-1} \\
t_2 & qt_1 & 2gq^2 & \ddots  & q^2t_{n-2} \\
\vdots & \ddots & \ddots & \ddots  & \vdots \\
t_n & qt_{n-1} & q^2t_{n-2} & \ldots  & 2gq^n \\
\end{pmatrix}.
\end{equation}
Being a Gram matrix, it must be positive semi-definite (PSD). 
Moreover, because the number of points $N_k$ of points of $X$ over $\F_{q^k}$ can only be greater or equal than $N_1$, the variables $t_k$ are also constrained, using \eqref{eq:deftk}, by the affine inequalities
\begin{equation}\label{eq:tk}
 t_k \leq t_1 + q^k -q.
\end{equation} 
Therefore, the vector $(t_1, \ldots, t_n)$ must belong to a spectrahedron, and by minimising $t_1$ over this convex set, one gets a lower bound on $t_1$ for any curve. 
In fact, this semi-definite program is closely related to the dual of the optimization problem solved by Oesterl\'e, as shown in \cite{HP19}, but provides a more geometric approach. 
Moreover, note that for $n=1$, one gets exactly the bound by Weil \eqref{eq:Weil1}, while $n=2$ leads to Ihara's
\cite{Ihara}. 
For this reason, for $n \geq 2$, these bounds can be seen as higher order Weil bounds, and we call this sequence of bounds the \textit{Weil-Oesterl\'e} hierarchy.
Nevertheless, if increasing $n$ leads in theory to better bounds, it was shown in \cite{HP19} that for any field size $q$, there is an explicit sequence $g_n = g_n(q)$ such that if $g_k \leq g \leq g_{k+1}$, then the bound for $q$ and $g$ does not improve for $k > n$.

While this approach hides the eigenvalues $\omega_j$ of the Frobenius, they play a crucial role in the next Section.

\subsection{Serre's improvement on Weil's bound}\label{sec:SerreWeil}

Turning back to Weil's bound \eqref{eq:Weil1}, a simple observation due to Serre leads to the following refinement. 
Recall that 
\[
t_1 = 1+q -N_1 = \sum_{j=1}^g \omega_j + \overline{\omega_j},
\]
with $\left| \omega_j \right| = \sqrt{q}$, and $\{ \omega_1, \ldots, \omega_g, \overline{\omega_1}, \ldots, \overline{\omega_g} \}$ stable under the action of $\Gal(\overline{\Q}/\Q)$.
If we set $\tau_1(\omega) = \omega + \overline{\omega}$, then for every $1 \leq j \leq g$, 
\[
\tau_1(\omega_j) = 2\sqrt{q}\cos(\theta_j) \in [-2\sqrt{q}, 2\sqrt{q}],
\]
so that 
\[
\tau_1(\omega_j) + \lfloor 2\sqrt{q} \rfloor + 1 >0,
\]
where $\lfloor 2\sqrt{q} \rfloor$ denotes the floor of $2\sqrt{q}$.  
Now, the arithmetic-geometric mean inequality implies that
\begin{equation}\label{eq:amgm}
\frac{1}{g}\sum_{j=1}^g (\tau_1(\omega_j) + \lfloor 2\sqrt{q} \rfloor + 1) \geq \left( \prod_{j=1}^g (\tau_1(\omega_j) + \lfloor 2\sqrt{q} \rfloor + 1) \right)^{\frac{1}{g}}.
\end{equation}
Next, observe that the product $\prod_{j=1}^g (\tau_1(\omega_j) + \lfloor 2\sqrt{q} \rfloor + 1)$ is 
\begin{enumerate}[label = \roman*)]
\item a rational number, because it is invariant under $\Gal(\overline{\Q}/\Q)$,
\item an algebraic integer, since it is a product of algebraic numbers,
\item positive, as a product of positive real numbers.
\end{enumerate}
This product thus has to be a positive integer, and is therefore at least $1$, so that \eqref{eq:amgm} gives
\begin{equation}\label{eq:SerreWeil}
t_1 \geq -g\lfloor 2\sqrt{q} \rfloor.
\end{equation}
If it does not improve upon Weil's bound when $q$ is a square, this refinement can lead to substantial improvements for non-square values of $q$ when $g$ grows.  
\section{Improvement on Ihara's bound}\label{sIS}

\subsection{A generalization of Serre's argument}

We first describe a general strategy to get an analogue of the trick described in Section~\ref{sec:SerreWeil} to higher order Weil-Oesterl\'e bounds. 
Recall that for any $k\geq 1$, 
\begin{equation}\label{eq:locglob}
t_k = \sum_{j=1}^g \tau_k(\omega_j)
\end{equation}
where 
\[
\tau_k(\omega_j) = \omega_j^k + \overline{\omega_j}^k = 2 q^\frac{k}{2} \cos(k \theta_j).
\]
Serre's argument relies on an affine inequality whose coefficients are integers, and which is satisfied by any $\tau_1(\omega) \in [-2\sqrt{q},2\sqrt{q}]$.
This idea can be generalized as follows.
\begin{lemma}\label{lem:GenSer}
Let $a_0, \ldots a_n$ integers. 
Assume that for every $\omega$ with $\left| \omega \right| = \sqrt{q}$,
\begin{equation}\label{eq:ineq}
\sum_{k=1}^n a_k\tau_k(\omega) + a_0> 0.
\end{equation}
Then 
\begin{equation}
\sum_{k=1}^n a_kt_k + ga_0 \geq g.
\end{equation}
\end{lemma}
\begin{proof}
Following the method given in Section~\ref{sec:SerreWeil}, we use the arithmetic-geometric mean, which gives
\begin{equation}\label{eq:genamgm}
\frac{1}{g}\sum_{j=1}^g \left(\sum_{k=1}^n a_k\tau_k(\omega_j) + a_0 \right) \geq \left( \prod_{j=1}^g \left(\sum_{k=1}^n a_k\tau_k(\omega_j)+ a_0 \right) \right)^{\frac{1}{g}}.
\end{equation}
Because the coefficients $a_k$ are integers, $\Gal(\overline{\Q}/\Q)$ permutes the algebraic integers $\sum_{k=1}^{\textcolor{violet}{n}} a_k\tau_k(\omega_j) + a_0$, for $1\leq j \leq g$.
Therefore the product of these positive real numbers is a positive integer, hence greater than $1$. 
Together with \eqref{eq:locglob}, this implies the result.
\end{proof}

When the coefficients $a_k$ are nonnegative for $1\leq k \leq n$, we can further combine Lemma~\ref{lem:GenSer} with the conditions \eqref{eq:tk}, to get the following general bound.

\begin{theorem}\label{thm:GenBound}
Let $a_0, \ldots a_n$ integers such that $a_k\geq 0$ for $1 \leq k \leq n$. 
Assume that for every $\omega$ with $\left| \omega \right| = \sqrt{q}$, the affine integral inequality
\begin{equation}\label{eq:integineq}
\sum_{k=1}^n a_k\tau_k(\omega) + a_0> 0
\end{equation}
holds.
Then 
\begin{equation}
t_1 \geq  \frac{g(1-a_0)-\sum_{k=1}^n a_k(q^k-q)}{\sum_{k=1}^n a_k}.
\end{equation}
\end{theorem}
\begin{proof}
The conditions of Lemma~\ref{lem:GenSer} are satisfied.
Then the conditions  $t_k \leq t_1 + q^k -q$ from \eqref{eq:tk} give, because $a_k \geq 0$ for $1\leq k \leq n$,
\[
\left(\sum_{k=1}^n a_k\right)t_1 + \sum_{k=1}^n a_k(q^k-q) \geq g(1-a_0), 
\]
which proves the result.
\end{proof}

In the following, we focus on the case $n=2$, explaining how to obtain inequalities of the form \eqref{eq:ineq} that improve upon Ihara's bound.

\subsection{Ihara's bound as the Weil-Oesterl\'e bound of order $2$}
\begin{figure}[H]
\begin{center}
\includegraphics[scale=1]{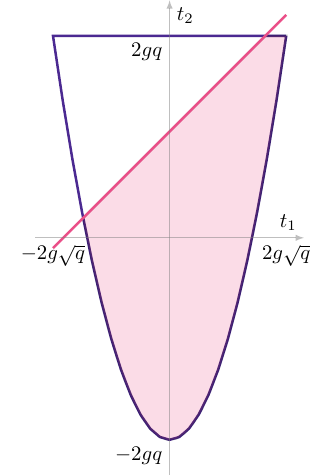}
\end{center}
\caption{The Weil domain for $n=2$.\label{fig:Ihara}}
\end{figure}

We first review how to obtain Ihara's bound with the approach from \cite{HP19}.
For $n=2$, the matrix in \eqref{eq:Toep} is 
\[
\Gram\left(p(\Gamma^0),p(\Gamma^1),p(\Gamma^2)\right)
=
\begin{pmatrix}
2g & t_1 & t_2  \\
t_1 & 2gq & qt_1  \\
t_2 & qt_1 & 2gq^2  \\
\end{pmatrix}.
\]
According to \cite{HP19}, a convenient change of basis make the computation easier:
\begin{equation}\label{eq:blocks}
\Gram\left(
\frac{q p(\Gamma^0) + p(\Gamma^2)}{\sqrt{2q}},
\frac{p(\Gamma^1)}{\sqrt{2q}},
\frac{q p(\Gamma^0) - p(\Gamma^2)}{\sqrt{2q}}
\right)
=
\begin{pmatrix}
2gq+t_2 & t_1 & 0  \\
t_1 & g & 0  \\
0 & 0 & 2gq-t_2  \\
\end{pmatrix}.
\end{equation}
This matrix is positive semi-definite if and only if the couple $(t_1,t_2)$ satisfies both the affine constraint $t_2 \leq 2qg$ and the quadratic constraint $t_2 \geq t_1^2/g - 2qg$, namely $(t_1,t_2)$ belongs to the convex set defined by the parabola and the horizontal line depicted in blue in Figure~\ref{fig:Ihara}.

The only additional affine constraint from \eqref{eq:tk} is $t_2 \leq t_1 + q^2-q$. 
When $g < \frac{\sqrt{q}(\sqrt{q}-1)}{2}$, the corresponding line does not meet this region, and one recovers the Weil-Oesterl\'e bound of first order, namely Weil's bound.
When $g \geq \frac{\sqrt{q}(\sqrt{q}-1)}{2}$, this line restricts the feasible domain to the convex set depicted in Figure~\ref{fig:Ihara}, that we call the \textit{Weil domain} of order $2$. 
The minimal~$t_1$ then occurs at the smallest intersection between this line and the parabola given by $t_2 = t_1^2/g - 2qg$. 
In other words, Ihara's bound is given by the smallest root of the polynomial $t^2/g - t - 2qg - q^2+q$ and therefore
\begin{equation}\label{eq:BoundIhara}
t_1 \geq g\left(\frac{1-\sqrt{1+8q+4(q^2-q)/g}}{2}\right).
\end{equation}
Furthermore, following \cite{HP19}, this bound does not get better for higher $n$ while $g \leq \frac{\sqrt{q}(q-1)}{\sqrt{2}}$.

\subsection{Getting affine integral inequalities}\label{sec:AffConstraints}

In order to apply Lemma~\ref{lem:GenSer}, we need affine inequalities  satisfied by every $\tau_1(\omega)$ and $\tau_2(\omega)$.
Here we provide a way to obtain such inequalities.
Denote by $M$ the first block of the Gram matrix \eqref{eq:blocks}.
It can be decomposed as
\begin{equation}\label{eq:DecompMomega}
M = \begin{pmatrix}
2gq+t_2 & t_1 \\
t_1 & g 
\end{pmatrix} = \sum_{j=1}^g M(\omega_j),
\end{equation}
where 
\begin{equation}\label{eq:Momega}
M(\omega) = \begin{pmatrix}
2q+\tau_2(\omega) & \tau_1(\omega) \\
\tau_1(\omega) & 1 
\end{pmatrix}.
\end{equation}
For every $\omega$ with $\left| \omega \right| = \sqrt{q}$, this is a PSD matrix of rank $1$, because
\[
\tau_1(\omega)^2 = (2\sqrt{q}\cos(\theta))^2 = 4q \frac{1+\cos(2\theta)}{2}  = 2q + \tau_2(\omega).
\] 
By duality, it follows that for every PSD matrix $A$, the trace inner product $\left\langle M(\omega), A \right\rangle = \Tr(M(\omega)A)$ is nonnegative for every $\omega$ with $\left| \omega \right| = \sqrt{q}$, and if moreover $A$ is definite positive, then the inequality becomes strict.
This provides a generic way to get inequalities that can be used in Lemma~\ref{lem:GenSer}.
\begin{lemma}\label{lem:ineqA2}
Let $A = \begin{pmatrix}
d & a \\ a & b
\end{pmatrix}$ be a positive definite matrix. Then for every $\omega$ with $\left| \omega \right| = \sqrt{q}$,
\[
d\tau_2(\omega) + 2a\tau_1(\omega) + 2qd + b >0.
\]
\end{lemma}
When $A$ runs through all possible positive definite matrices, Lemma~\ref{lem:GenSer} gives affine equations that have to be satisfied by $(t_1,t_2)$.

Then, going from Lemma~\ref{lem:GenSer} to Theorem~\ref{thm:GenBound} corresponds to taking into consideration the additional constraint $t_2 \leq t_1 + q^2 -q$, which gives the following.
\begin{theorem}\label{thm:boundA2}
Let $A = \begin{pmatrix}
d & a \\ a & b
\end{pmatrix}$ be a positive definite matrix. Assume that $d, b\in \N$ and $2a\in \N$.
Then 
\[
t_1 \geq \frac{g(1-2qd-b)-d(q^2-q)}{d+2a}.
\]
\end{theorem}

\begin{figure}[H]
\begin{center}
\includegraphics[scale=.9]{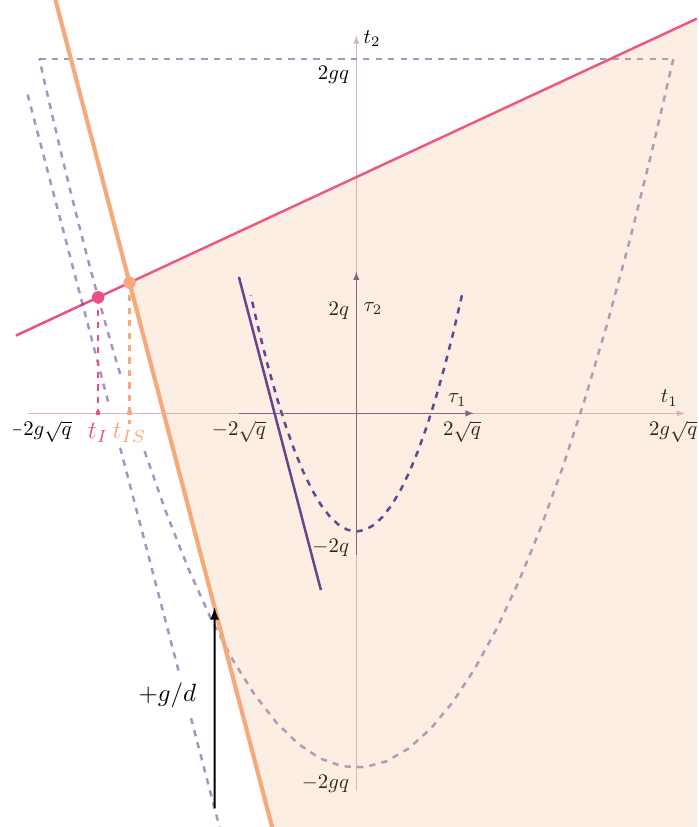}
\end{center}
\caption{For every $j$, the couple $(\tau_1(\omega_j),\tau_2(\omega_j))$ lies on the interior parabola.
Every matrix $A$ satisfying the assumptions of Theorem~\ref{thm:boundA2} first gives an affine constraint on $(\tau_1(\omega_j),\tau_2(\omega_j))$, such as the one represented by the plain interior blue line. 
By summing over the $\omega_j$'s, the point $(t_1,t_2)$ is thus constrained by the dashed exterior blue line.
However, the Serre type Lemma~\ref{lem:GenSer} ensures that $(t_1,t_2)$ is in fact above the shifted orange line.
If we can find a matrix $A$ such that this line intersects the magenta line given by the additional constraint $t_2 \leq t_1 + q^2 -q$ in the interior of the Weil domain, we get a better lower bound on $t_1$.
 \label{fig:IharaTrunc}}
\end{figure}

Figure~\ref{fig:IharaTrunc} sums up the situation so far: from every matrix $A$ that fulfills the conditions of Theorem~\ref{thm:boundA2} we obtain an affine constraint on $(t_1,t_2)$, that might exclude some region from the Weil domain.

\subsection{Optimisation of the bound}

It remains to prove that we can find a matrix $A$ such that Theorem~\ref{thm:boundA2} gives a better bound than Ihara's bound. 
The next theorem, main result of the paper, shows that there is a choice of $A$ that improves upon Ihara's bound, and evaluates the gain between the two bounds.

\begin{theorem}\label{thm:main}
Let~$X$ be an an absolutely, irreducible, smooth, projective curve, of genus~$g$, defined over the finite field~$\F_q$.
Suppose that~$g \geq \frac{\sqrt{q}(\sqrt{q}-1)}{2}$, let~$t_I$ be the Ihara trace
\begin{align*}
&t_I = g\left(\frac{1-\sqrt{1+8q+4(q^2-q)/g}}{2}\right),
&
&\text{and set:}
&
&\alpha = - \frac{t_I}{g}.
\end{align*}
Then, the best upper bound obtained using theorem~\ref{thm:boundA2} with a matrix having an upper left coefficient equal to~$1$ is
for
\begin{align*}
&A =
\begin{pmatrix}
1 & a \\ a & \lfloor a^2 \rfloor + 1
\end{pmatrix},
&
&\text{where}
&
&a = \lfloor \alpha \rfloor +  \frac{1}{2} = \left\lfloor - \frac{t_I}{g} \right\rfloor +  \frac{1}{2}.
\end{align*}
The Ihara's upper bound improved {\it à la Serre} is
$$
\sharp X(\F_q)
\leq (q+1) - t_I - g\frac{(\alpha - \lfloor \alpha \rfloor)(\lceil \alpha \rceil -\alpha)}{2\lceil \alpha \rceil}
$$
and the gain between the improved bound and the original one equals
\begin{equation}\label{eq:diffmain}
g\frac{(\alpha - \lfloor \alpha \rfloor)(\lceil \alpha \rceil -\alpha)}{2\lceil \alpha \rceil}.
\end{equation}
\end{theorem}

\begin{proof}

First note for a matrix of the form $A = \begin{pmatrix}
1 & a \\ a & b
\end{pmatrix}$, the lower bound given by Theorem~\ref{thm:boundA2} gets only worse when $b$ increases, therefore it is better to take $b$ as the smallest integer such that $A$ is positive definite, namely
$$
b = \lfloor a^2 \rfloor + 1.
$$
Let us denote by~$\newlowerbound(a)$, where $IS$ stands for Ihara-Serre, the lower bound on~$t_1$ obtained using Theorem~\ref{thm:boundA2} with~$d=1$,
$a$ such that~$2a\in\N$, and~$b = \lfloor a^2 \rfloor + 1$.
According to Theorem~\ref{thm:boundA2}, 
\[
\newlowerbound(a) = -g\left(\frac{\lfloor a^2 \rfloor + 2q + (q^2-q)/g}{1+2a}\right),
\]
and therefore
\[
\frac{\newlowerbound(a)-t_I}{g}= - \frac{\lfloor a^2 \rfloor + 2q + (q^2-q)/g}{1+2a} + \alpha = - \frac{\lfloor a^2 \rfloor + 2q + (q^2-q)/g - \alpha - 2 a \alpha}{1+2a}.
\]
We observe that since 
\[\alpha = \frac{\sqrt{1+8q+4(q^2-q)/g}-1}{2},\]
we can rewrite
\[
2q + (q^2-q)/g = \alpha(\alpha + 1)
\]
and thus
$$
\newlowerbound(a)-t_I=\frac{-g}{1+2a}\left(\alpha^2-2a\alpha+\lfloor a^2 \rfloor\right).
$$
Since $2a \in \N$, we have two distinct cases.

$\bullet$ If $2a$ is even, then $a\in \N$, and thus $\lfloor a^2 \rfloor = a^2$, and in this case 
\[
\newlowerbound(a)-t_I = \frac{-g}{1+2a}\left(\alpha - a \right)^2 <0,
\]
we always get a weaker bound.

$\bullet$ If $2a$ is odd, then we can write $a = k +1/2$ with $k\in\N$, and in this case $\lfloor a^2 \rfloor = k(k+1)$. 
Then the bounds reads
\begin{align*}
\newlowerbound(a)-t_I &= \frac{-g}{2(1+k)}\left(\alpha^2 - (2k+1)\alpha + k(k+1) \right)
\\  &= \frac{-g}{2(1+k)}\left(\alpha - k \right)\left(\alpha - k - 1 \right)
\\  &= \frac{g}{2(1+k)}\left(\alpha - k \right)\left(k+1 - \alpha \right)
,
\end{align*}
and this difference is positive if and only if $k < \alpha < k+1$.
If $\alpha$ is not an integer, this occurs only for $k = \lfloor \alpha \rfloor$.
In the very specific case where $\alpha$ is an integer, then we recover Ihara's bound for $k = \alpha$ and $k = \alpha -1$.
\end{proof}

\subsection{Improvements on previous bounds}\label{sEtudeGain}
Theorem~\ref{thm:main} gives an explicit comparison between our bound and Ihara's bound.
In this section we comment on this comparison, with several points of view.

\subsubsection{Fixed $q$, $g$ grows}\label{sec:qfix}
First, observe that for a fixed $q$, when $g$ grows, 
\[
\alpha = -\frac{t_I}{g} = \frac{\sqrt{1+8q+4(q^2-q)/g}-1}{2}
\]
goes to 
\[
\alpha^\infty = \frac{\sqrt{1+8q}-1}{2}.
\]
Thus, when $g$ goes to infinity, our Ihara-Serre bound 
\begin{equation}
\newlowerbound=t_I + g\frac{(\alpha - \lfloor \alpha \rfloor)(\lceil \alpha \rceil -\alpha)}{2\lceil \alpha \rceil}
\end{equation}
on $t_1$ satisfies
\[
\frac{\newlowerbound-t_I}{g} \to_{g\to \infty} \frac{(\alpha^\infty - \lfloor \alpha^\infty \rfloor)(\lceil \alpha^\infty \rceil -\alpha^\infty)}{2\lceil \alpha^\infty \rceil}.
\]
By pushing the analysis further, one finds that our gain $\newlowerbound-t_I$, as a function of $g$, has an asymptote when $g$ goes to infinity, whose constant term is $\displaystyle{\frac{(\lfloor \alpha^\infty \rfloor  - \alpha^\infty +1/2)(q^2-q)}{(\lfloor \alpha^\infty \rfloor + 1)\sqrt{8q+1}}}$. 
This has two consequences. 
First it shows that whenever $\alpha^\infty$ is not an integer (namely whenever $q \neq 3$), the difference between our bound and Ihara's bound goes to infinity. 
Second, this shows that our method gives an improvement on the upper bound given by Ihara's bound on the constant 
\[A(q) = \limsup_{g\to\infty} \frac{N_1 - (q+1)}{g}. \]
When Ihara's bound gives $A(q) \leq \alpha^\infty$, our bound gives
\begin{equation}
A(q) \leq \alpha^\infty - \frac{(\alpha^\infty - \lfloor \alpha^\infty \rfloor)(\lceil \alpha^\infty \rceil -\alpha^\infty)}{2\lceil \alpha^\infty \rceil}.
\end{equation}

Of course our bound, as Ihara's, gets weaker than higher order Weil-Oesterl\'e bounds when $g$ grows, and it was shown in \cite{HP19} that these bounds recover the Drinfeld-Vl{\u a}du\c t \cite{VD} bound  when $n$ grows, namely $A(q) \leq \sqrt{q} - 1$.
However, it makes sense to compare our bound with Ihara's, since they correspond to Weil-Oesterl\'e bounds of the same order.

\subsubsection{The Ihara range}\label{sec:CompIhara} For every $q$ there is a range $[g_2,g_3]$ for which Ihara's bound is not improved by higher order Weil-Oesterl\'e bounds.
We now focus on this range, in order to show that our bounds gives explicit improvements on Ihara's bound, even when it was the best bound known so far.
Recall that $g_2 = \left\lceil \frac{\sqrt{q}(\sqrt{q}-1)}{2} \right\rceil$ and $g_3 = \left\lfloor \frac{\sqrt{q}(q-1)}{\sqrt{2}} \right\rfloor$.

\begin{figure}[!h]
\begin{tabular}{cc}
\includegraphics[scale=.38]{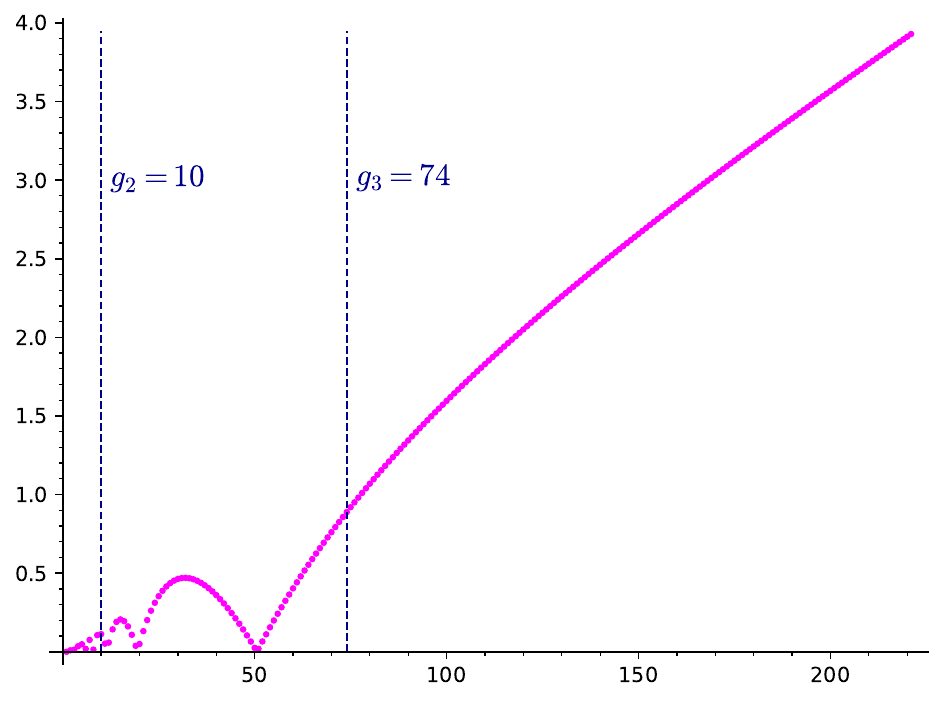}
&\includegraphics[scale=.38]{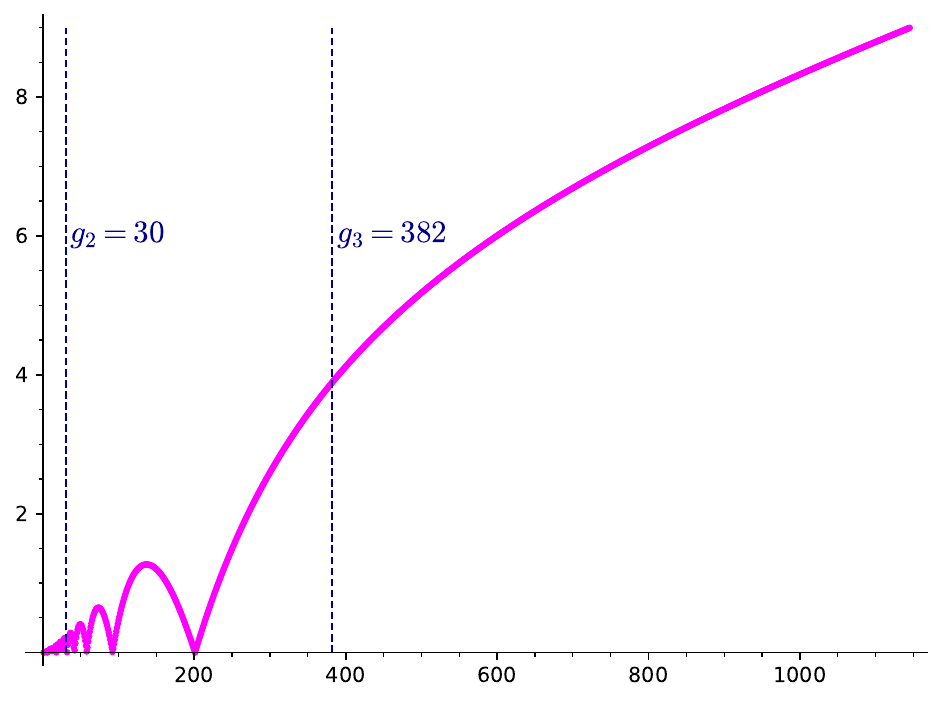}
\\
$ q = 23$
&
 $q= 67$
\\\includegraphics[scale=.38]{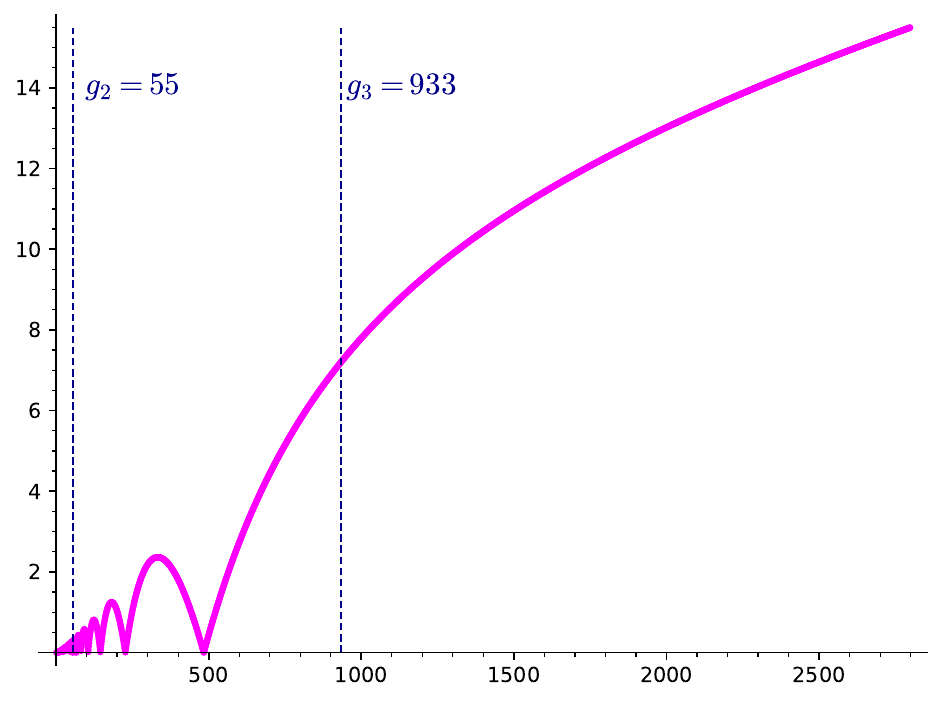}
&\includegraphics[scale=.38]{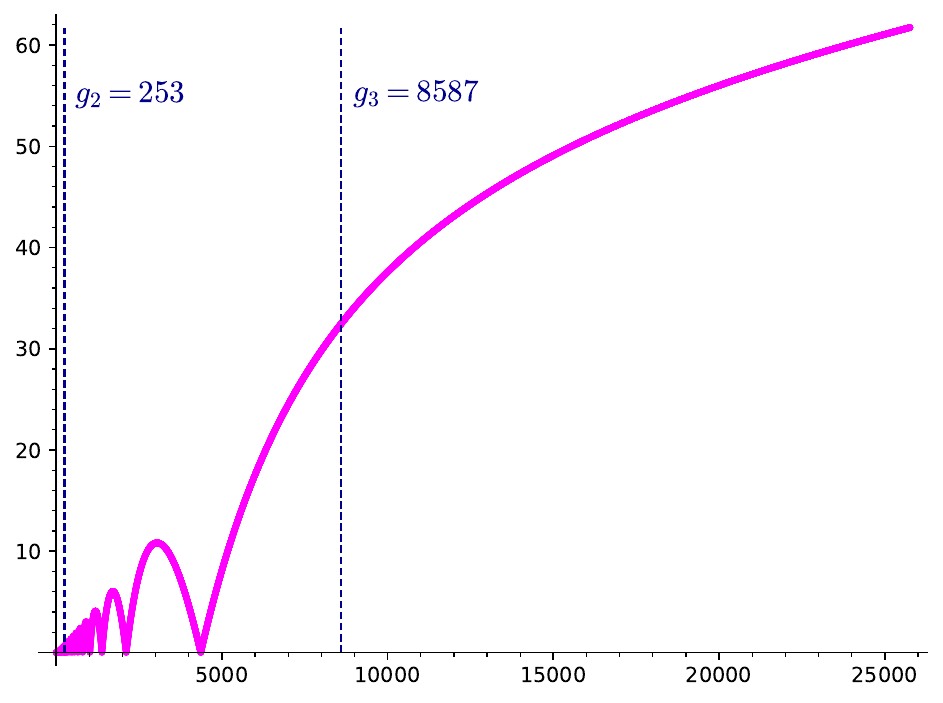}
\\
 $q = 121=11^2$

&
$q = 529=23^2$
\end{tabular}
\caption{The difference \eqref{eq:diffmain} for some values of $q$ and $g \leq 3g_3$. \label{fig:plots}}
\end{figure}

In Figure~\ref{fig:plots}, we plot the gain \eqref{eq:diffmain} of our bound compared with Ihara's bound for several values of $q$, and in each case, values of $g$ up to $3g_3$.
One can observe that in comparison with the asymptotic behavior, the ceiling and floor functions in \eqref{eq:diffmain} make the the analysis of the gain difficult in the range $[g_2,g_3]$, even though the difference seems to increase with~$q$.
Furthermore, note that our bound improves upon Ihara's for several values of~$g$ even when~$q$ is a square, while Serre's trick does not improve upon Weil's bound in that case.

\subsubsection{A sequence $g_q$ with $q$ increasing}\label{sec:CompSeq}
\begin{figure}[!h]
\includegraphics[scale=.6]{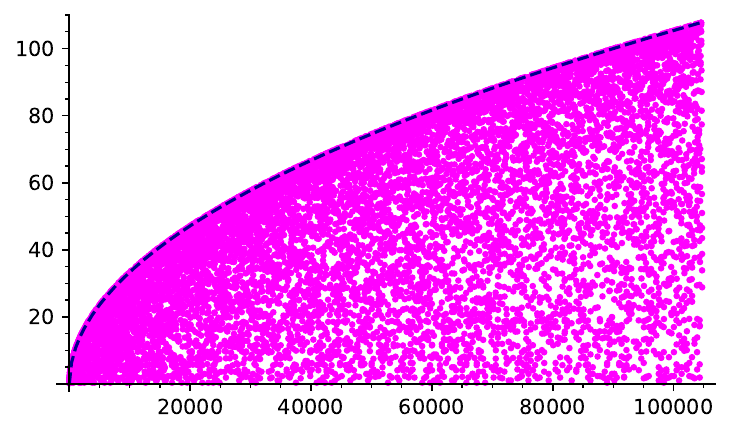}
\caption{The difference \eqref{eq:explicitqg} for $q$ prime, and $g = 4q$. For comparison, the dashed curve is the function $\sqrt{q}/3$. \label{fig:stats}}
\end{figure}

To make this more explicit, we study an explicit sequence where the genus  $g_q$ is in the Ihara range $[g_2,g_3]$ for $q$, and $q$ increases.
More precisely, let for instance $g_q = 4q$.
Then $g_2 \leq g_q \leq g_3$ whenever $q \geq 34$, and $\alpha$ in Theorem~\ref{thm:main} simplifies to
\[
\alpha = \frac{3\sqrt{q}-1}{2},
\]
and \eqref{eq:diffmain} becomes
\begin{equation}\label{eq:explicitqg}
\frac{2q}{\lceil \frac{3\sqrt{q}-1}{2} \rceil}(\alpha - \lfloor \alpha \rfloor)(\lceil \alpha \rceil -\alpha).
\end{equation}
Moreover, since the maximum of the map $x\mapsto (x - \lfloor x \rfloor)(\lceil x \rceil -x)$ is $1/4$ when $x$ is a half-integer, the best gain one can hope for $g_q$ is  
\[
\frac{q}{2\lceil \frac{3\sqrt{q}-1}{2} \rceil}.
\]
In Figure~\ref{fig:stats} we plot, for the $10000$ first prime numbers, the gain achieved with our bound compared with Ihara bound. 
One can then see that for a high proportion of~$q$, the improvement is close to $\sqrt{q}/3$ when $q$ grows.

Also, if we take $q = 2^{2k}$ for $k$ a large enough integer, then 

\[
\alpha = \frac{3 \cdot 2^k - 1}{2} = 3\cdot 2^{k-1} - \frac{1}{2}
\]

and in this case the gain is optimal, namely $\sqrt{q}/3$.
This shows in particular that our bound improves upon Ihara's for infinitely many couples $(q,g)$ where $g$ is in the Ihara range.

\subsubsection{Small $q$: comparison with manYPoints} \label{sec:MP}
One can see in Figure~\ref{fig:plots} that for small values of $q$, the improvement in the Ihara range can be small.
However, since a lower bound on $t_1$ gives an upper bound on the number of points $N_1$ of a curve $X$ over $\F_q$, if the bound obtained by Ihara is close to its floor, our bound can sometimes go below this floor and therefore provide an improvement by $1$ in the upper bound.
In Table~\ref{tab:MP} we give the pairs $(g,q)$ for which our bound improves upon Ihara's bound in this sense, among the values of $q$ and $g$ displayed in manYPoints \cite{MP}.
However, for such small values of $g$ and $q$, other specific methods were developed \cite{HL03, HL12}, and give sometimes better bounds than Ihara's, even in the Ihara range.
If our bound often meets these improvements with our generic method, in some cases our bound is worse.
Also, for $g$ close to $g_2$, Serre's improvement on Weil's bound can be much stronger than Ihara's bound, and also beat our bound.
These cases where our bound does not meet the current records are displayed in parentheses. 
Nevertheless, our result provides several new records: for the numbers displayed in bold in Table~\ref{tab:MP}, our upper bound on $N_1$ improves upon the current records by $1$.

\begin{center}
\begin{table}[!h]
\begin{tabular}{ c    c   }
$q$  &  $g$  \\ \hline 
\rule{0pt}{4mm}%
$9$ & $(12)$,  $(17)$\\
$11$ & $8$,  $19$, $23$, $24$\\
$13$ & $16$,  $19$, $(22)$, $(25)$\\
$16$ & $13$,  $20$, $34$, $35$, $39$, $40$, $41$\\
$17$ & $15$,  $17$, $22$, $24$, $29$, $42$, $45$\\
$19$ & $12$,  $23$, $27$, $30$, $31$, $34$, $35$, $37$, $38$, $41$, $42$, $44$, $45$, $48$, $49$\\
$23$ & $14$,  $22$, $27$, $30$, $32$, $35$, $38$, $43$, $46$\\
$25$ & $20$,  $31$, $34$, $37$, $39$, $40$, $42$, $43$, $45$, $47$, $48$, $50$\\
$27$ & $(11)$,  $14$, $48$, $49$, $50$\\
$29$ & $18$,  $30$, $34$, $35$, $38$, $39$, $43$, $44$, $48$\\
$31$ & $24$,  $27$, $29$, $41$, $43$, $45$, $47$, $49$, $50$\\
$32$ & $26$,  $29$, $41$, $46$, $48$, $50$\\
$37$ & $28$,  $31$, $34$, $41$, $45$, $46$, $49$, $50$\\
$41$ & $(18)$,  $29$, $30$, $39$, $40$, $43$, $44$, $47$, $50$\\
$43$ & $(21)$,  $30$, $31$, $32$, $46$, $(47)$, $48$, $49$, $50$\\
$47$ & $(23)$,  $32$, $38$, $40$, $42$, $45$, $47$, $\textbf{50}$\\
$49$ & $33$,  $37$, $46$, $\textbf{49}$\\
$53$ & $33$,  $38$, $\textbf{47}$, $\textbf{48}$, $\textbf{49}$, $\textbf{50}$\\
$59$ & $32$,  $33$, $34$,  $40$,  $43$, $44$, $\textbf{47}$, $\textbf{50}$\\
$61$ & $(27)$,  $\textbf{48}$, $\textbf{49}$, $\textbf{50}$\\
$64$ & $39$,  $43$, $47$\\
$67$ & $36$,  $\textbf{44}$, $\textbf{46}$, $\textbf{48}$, $\textbf{50}$\\
$71$ & $(32)$,  $41$\\
$73$ & $35$,  $38$, $43$\\
$79$ & $(38)$,  $(39)$, $\textbf{49}$\\
$81$ & $\textbf{50}$\\
$83$ & $\textbf{43}$,  $\textbf{44}$\\
$89$ & $(42)$,  $(45)$, $\textbf{50}$\\
$97$ & $(46)$,  $(49)$\\
\hline
\end{tabular}
\caption{The couples $(q,g)$ with $q \leq 100$, $g\leq 50$ such that our bound on $N_1$ is strictly better than Ihara's. In bold, the new records compared with \cite{MP}. In parentheses, better bounds can be found, either by Serre's improvement on Weil's bound, or by the techniques from \cite{HL03, HL12}.
For all the other numbers, our bound meets the current record.}
\label{tab:MP}
\end{table}
\end{center}

\section{Extension to higher order Weil-Oesterl\'e bounds}\label{sec:extensions}

In principle, our approach can be generalized to Weil-Oesterl\'e bounds of higher order. 
However, it becomes harder to obtain a statement similar to Theorem~\ref{thm:main}, where we can explicitely optimise over the matrix $A$ and give a closed formula for the gain compared to the corresponding Weil-Oesterl\'e bound.
Another, more experimental, approach consists, for a fixed couple $(q,g)$, in trying several matrices $A$ and search for the one that provides the best bound.

\subsection{Experimental search}
We sketch this approach for the Weil-Oesterl\'e bound of order $3$, and show that it can successfully improve upon Oesterl\'e.
For $n=3$, the Gram matrix in \eqref{eq:Toep} is 
\[
\Gram\left(p(\Gamma^0),p(\Gamma^1),p(\Gamma^2),p(\Gamma^3)\right)
=
\begin{pmatrix}
2g & t_1 & t_2 & t_3 \\
t_1 & 2gq & qt_1 & qt_2 \\
t_2 & qt_1 & 2gq^2 & q^2t_1 \\
t_3 & qt_2 & q^2t_1 & 2gq^3 \\
\end{pmatrix}.
\]
As for Ihara, after a convenient change of basis, the Gram matrix 
\begin{align*}
\Gram
&\left(
\frac{q^{\frac{3}{2}} p(\Gamma^0) + p(\Gamma^3)}{\sqrt{2}q^{\frac{3}{4}}},
\frac{q^{\frac{1}{2}} p(\Gamma^1) + p(\Gamma^2)}{\sqrt{2}q^{\frac{3}{4}}},
\frac{q^{\frac{1}{2}} p(\Gamma^1) - p(\Gamma^2)}{\sqrt{2}q^{\frac{3}{4}}},
\frac{q^{\frac{3}{2}} p(\Gamma^0) - p(\Gamma^3)}{\sqrt{2}q^{\frac{3}{4}}}
\right)\\
&=
\begin{pmatrix}
2gq^\frac{3}{2}+t_3 & \sqrt{q}t_1 + t_2 & 0 & 0  \\
\sqrt{q}t_1 + t_2 & 2g\sqrt{q}+t_1 & 0 & 0  \\
0 & 0 & 2g\sqrt{q}-t_1 & \sqrt{q}t_1 + t_2    \\
0 & 0 & \sqrt{q}t_1 + t_2 & 2gq^\frac{3}{2}-t_3   \\
\end{pmatrix} 
\end{align*}
is of course a positive semidefinite matrix, and makes the computations easier.
Let 
\[
M = \begin{pmatrix}
2gq^\frac{3}{2}+t_3 & \sqrt{q}t_1 + t_2   \\
\sqrt{q}t_1 + t_2 & 2g\sqrt{q}+t_1  \\
\end{pmatrix}
\]
be the first block of this matrix.
Then, the Weil-Oesterl\'e bound of order $3$, optimal for $g$ in the range $[g_3,g_4]$, is given by the intersection with minimal $t_1$ between the hypersurface given by $\det(M) = 0$ and the line given by the equations $t_2 = t_1 + q^2- q$ and $t_3 = t_1 + q^3-q$, see \cite{HP19}.  

We follow the strategy of Section~\ref{sec:AffConstraints}, and write $M = \sum_{j=1}^g M(\omega_j)$, with 
\[
M(\omega) = \begin{pmatrix}
2q^\frac{3}{2}+\tau_3(\omega) & \sqrt{q}\tau_1(\omega) + \tau_2(\omega)   \\
\sqrt{q}\tau_1(\omega) + \tau_2(\omega) & 2\sqrt{q}+\tau_1(\omega)  \\
\end{pmatrix}.
\]
Again, this matrix is a PSD matrix of rank $1$, because in addition to $\tau_2(\omega)= \tau_1(\omega)^2 - 2q$ we have $\tau_3(\omega)= \tau_1(\omega)^3 - 3q\tau_1(\omega)$.
In this context, we can adapt Lemma~\ref{lem:ineqA2} and Theorem~\ref{thm:boundA2} to get the following statement.

\begin{theorem}\label{thm:boundA3}
Let $A = \begin{pmatrix}
d & a \\ a & b
\end{pmatrix}$ be a positive semi-definite matrix. Assume that $d, 2a$, and $2a\sqrt{q} + b$ are natural integers.
Then 
\[
t_1 \geq \frac{-g \left\lfloor  2q^{3/2}d +2b\sqrt{q} \right\rfloor-d(q^3-q)-2a(q^2-q)}{d+2a + 2a\sqrt{q}+b}. \]
\end{theorem}
\begin{proof}
Because $A = \begin{pmatrix}
d & a \\ a & b
\end{pmatrix}$ and $M$ are both PSD, we get for every $\omega$ with $\left| \omega \right| = \sqrt{q}$ the inequality
\[
 d \tau_3(\omega) + 2a \tau_2(\omega) + (2a\sqrt{q} + b)\tau_1(\omega) +  2q^{3/2}d +2b\sqrt{q} \geq 0,
 \]
which implies the strict inequality
\[
 d \tau_3(\omega) + 2a \tau_2(\omega) + (2a\sqrt{q} + b)\tau_1(\omega) + \left\lfloor  2q^{3/2}d +2b\sqrt{q} \right\rfloor + 1 > 0. \]
The assumptions ensuring that all the coefficients are integers, we can apply Lemma~\ref{lem:GenSer}, which yields
\[
 d t_3 + 2a t_2 + (2a\sqrt{q} + b)t_1  \geq -g \left\lfloor  2q^{3/2}d +2b\sqrt{q} \right\rfloor. \]
As usual, since the coefficients are assumed to be nonnegative we conclude by adding the constraints $t_2 \leq t_1 + q^2 - q$ and $t_3 \leq t_1 + q^3 - q$.
\end{proof}

Using a computer algebra system, one can then compute the bound given by Theorem~\ref{thm:boundA3} for a large number of $A$, and check whether it is better than Oesterl\'e's bound. 
We can find such matrices $A$'s, and sometimes this leads, like in Section~\ref{sec:MP}, to a better bound on the number of points $N_1$ of a curve of genus $g$ over~$\F_q$. 
This even provides new bounds with respect to \cite{MP}.

\begin{theorem}\label{thm:rec3}
The number~$N_q(g)$ is bounded above by 
\[
N_q(g) \leq \begin{cases}
53 & \text{ if } q = 5 , g= 19, \\ 
76 & \text{ if } q = 7, g= 21, \\
129 & \text{ if } q = 8, g= 36, \\
163 & \text{ if } q = 11, g= 35. \\
\end{cases}
\]
\end{theorem}
\begin{proof} We apply Theorem~\ref{thm:boundA3} with the matrices $A$ given in the following table.

\begin{tabular}{ c  c  c  c   }
$(q,g)$  &  $A$ & lower bound on $t_1$ & upper bound on $N_1$  \\ \hline 
\rule{0pt}{4mm}%
$(5,19)$ & $\begin{pmatrix}
1 & 7/2 \\ 7/2 & -7\sqrt{5}+28
\end{pmatrix}$ &  $-\frac{1723}{36}$ &  $\frac{1939}{36} < 54$ \\
$(7,21)$ & $\begin{pmatrix}
3 & 29/2 \\ 29/2 & -29\sqrt{7}+147
\end{pmatrix}$ &  $-\frac{12348}{179}$ &  $\frac{13780}{179} < 77$ \\
$(8,36)$ & $\begin{pmatrix}
3 & 15 \\ 15 & -60\sqrt{2}+160
\end{pmatrix}$ &  $-\frac{23352}{193}$ &  $\frac{25089}{193} < 130$ \\
$(11,35)$ & $\begin{pmatrix}
2 & 14 \\ 14 & -28\sqrt{11}+191
\end{pmatrix}$ &  $-\frac{33580}{221}$ &  $\frac{36232}{221} < 164$ \\
\hline
\end{tabular}

\end{proof}

\subsection{Obstacles and limits}
There are several reasons that make the task of finding a closed formula as in Theorem~\ref{thm:main} hard for $n\geq 3$. 

On the positive side, we have an intuition on how to reduce the domain in which we search for a good matrix~$A$.
Indeed, if we denote by $t_W = (t_1,t_2,t_3)$ the solution of the initial Weil-Oesterl\'e problem, then the matrix $M$ has rank $1$ at $t_W$. 
If $v$ is a kernel vector of $M$, then the PSD matrix $A_0 = v^t v$ satisfies $\langle M, A_0 \rangle = 0$.
Geometrically, this equation defines the tangent space of the hypersurface defined by $\det(M)=0$ at~$t_W$. 
One can then look for a matrix $A$ close to $A_0$ which satisfies the conditions of Lemma~\ref{thm:boundA3}, this is how we obtained the matrices in Theorem~\ref{thm:rec3}. 

However the vector $v$, and therefore the matrix $A_0$, is unique only up to rescaling, and while fixing the top left coefficient of $A$ to $1$ appears to be optimal for $n=2$, this does not seem to be the case for $n =3$ (see the matrices in the proof of Theorem~\ref{thm:rec3}).
Furthermore, the way to approximate $A_0$, even after fixing the top left coefficient of $A$, seems to depend on the arithmetic of the numbers involved when taking floors.
This prevents us to find a canonical candidate for the optimal~$A$ as we were able to do for $n=2$.
More generally, when $n$ grows, the number of coefficients in $A$ increases, the relations between these coefficients become more complicated, and the initial optimal point $t_W$ is defined by a polynomial equation of increasing degree. 
All these reasons make the generalisation of Theorem~\ref{thm:main} a problem that seems difficult.

On the other hand, the numerical experiments we made, even for larger $q$ and~$g$, suggest that the improvement compared to the corresponding standard Weil-Oesterl\'e bound is weaker when $n=3$ compared to $n=2$, which is itself weaker than the improvement by Serre upon Weil's bound.
This might be explained by the fact that when $n$ grows, the spectrahedron
defined by the Gram matrix of order $n$ is described by equations of higher degree in more variables, and the truncation of this spectrahedron that we obtain with our affine inequalities becomes less powerful.

\bibliographystyle{amsalpha}
\bibliography{IharaSerre}
\end{document}